\documentclass[12pt,reqno]{amsart}
\usepackage{fullpage,xypic,xcolor,amssymb,amsthm}
\usepackage{mathpazo,verbatim,graphicx,float}
\usepackage[figurename=Diagram]{caption}
\usepackage[T1]{fontenc}
\newcommand{\Aut}{\text{Aut}}
\newcommand{\Out}{\text{Out}}
\newcommand{\Inn}{\text{Inn}}
\renewcommand{\P}{\mathbb P}
\newcommand{\complex}{\mathbb C}
\newcommand{\Hex}{\mathcal H}
\newcommand{\ZZ}{\mathbb Z}
\newcommand{\FF}{\mathbb F}
\newcommand{\frakp}{\mathfrak p} 
\newcommand{\Planedual}{({\mathbb P}^2)^*}
\newcommand{\demo}{{\sc Proof. }}
\newcommand{\lra}{\longrightarrow} 
\newcommand{\ra}{\rightarrow} 
\newcommand{\conic}{\mathbb K}
\newcommand{\SG}{\mathbb S}
\newcommand{\IN}{\mathcal N} 
\newcommand{\six}{\textsc{six}}
\newcommand{\ltr}{\textsc{ltr}}
\newcommand{\pasc}[6]{\left\{\begin{array}{ccc} #1 & #2 & #3\\ #4 &
     #5 & #6 \end{array} \right\}}
\newcommand{\transp}[2]{(#1 \, #2)}
\newcommand{\cyclett}[6]{(#1 \, #2) \, (#3 \, #4) \, (#5 \, #6)}
\newcommand{\synth}[6]{#1  #2. #3 #4. #5 #6}
\newcommand{\bbA}{\mathbb{A}}
\newcommand{\bbB}{\mathbb{B}}
\newcommand{\bbC}{\mathbb{C}}
\newcommand{\bbD}{\mathbb{D}}
\newcommand{\bbE}{\mathbb{E}}
\newcommand{\bbF}{\mathbb{F}}
\newtheorem{Theorem}{Theorem}[section]

\newtheorem{Proposition}[Theorem]{Proposition}
\begin{document} 
\title{On the enumerative geometry of Pascal's hexagram}
\author{Jaydeep Chipalkatti$^*$}
\thanks{*Department of Mathematics, University of Manitoba, Winnipeg,
  MB R3T 2N2. Canada. \\ Email: jaydeep.chipalkatti@umanitoba.ca}
\maketitle

\parbox{17cm}{ \small
{\sc Abstract:} 
Given six points $A,B,C,D,E,F$ on a nonsingular conic in the
complex projective plane, Pascal's theorem says that the three intersection points $AE \cap BF,
BD \cap CE, AD \cap CF$ are collinear. The line containing them is
called a pascal, and we get altogether $60$ such lines by permuting
the points. In this paper, we consider the enumerative problem of
finding the number of sextuples $(A, B, \dots, F)$ which correspond to
three pre-specified pascals. We use computational techniques in
commutative algebra to solve this problem in all cases. The results are
tabulated using the so-called 'dual' notation for pascals, which
is based upon the outer automorphism of $\SG_6$.} 

\bigskip 

Keywords: Pascal lines, Steiner's and Kirkman's theorems, intersection numbers. 

\bigskip 

AMS subject classification (2020): 14N05, 14N10. 

\bigskip 

\tableofcontents

\bigskip 

\section{Introduction} 
The object of this paper is to consider a certain kind of enumerative
problem which arises naturally in the context of Pascal's theorem. We
begin with an elementary introduction to the subject; the results are
described in Section~\ref{section.enumerative} once the required
notation is available. 
\subsection{Pascal's Theorem}
\begin{figure}
\includegraphics[width=9cm]{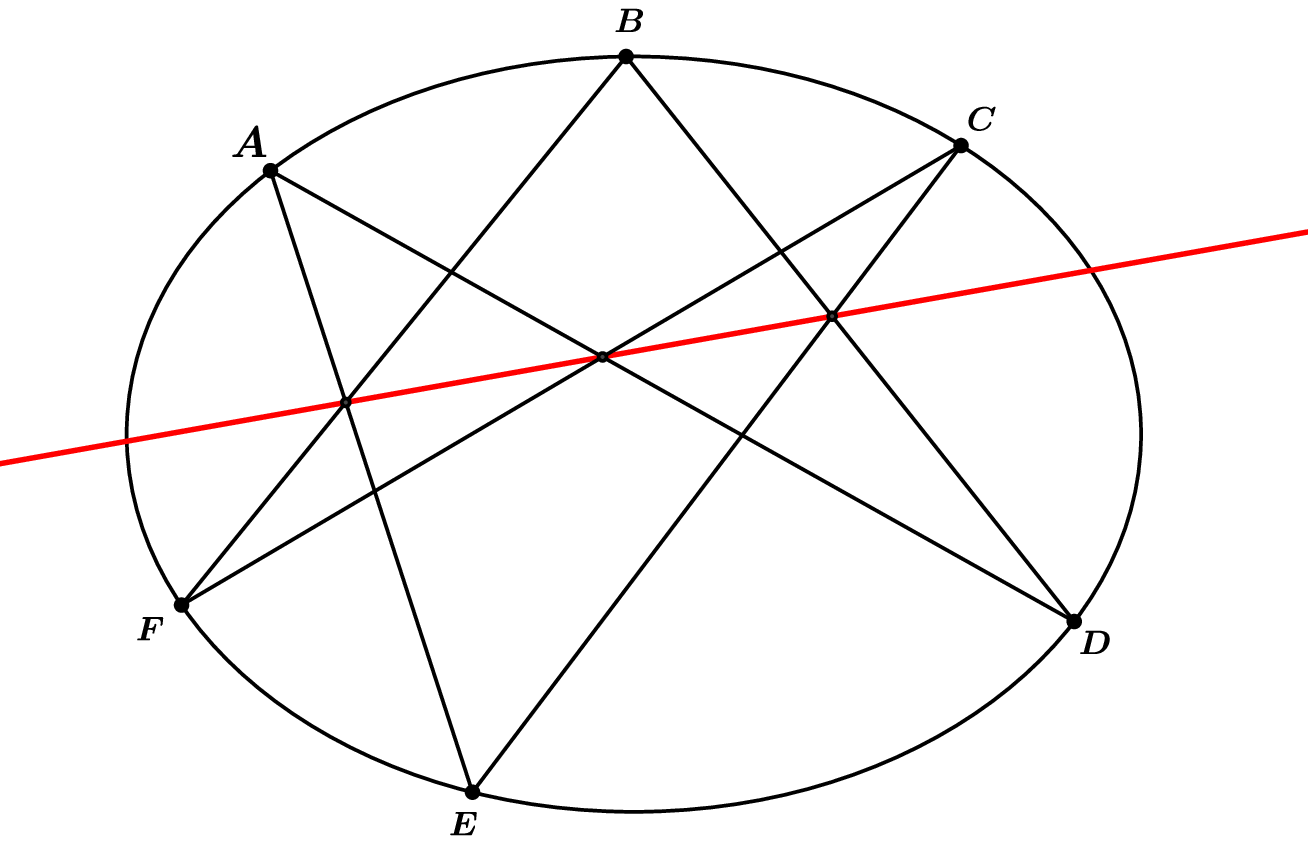}
\caption{Pascal's theorem} 
\label{diag:pascal} 
\end{figure}

Fix a nonsingular conic $\conic$ in the complex projective
plane $\P^2$. Suppose that we are given six distinct points
$A,B,C,D,E,F$ on $\conic$, 
arranged into an array
$\left[ \begin{array}{ccc} A & B & C \\ F & E & D \end{array}
\right]$. Then Pascal's theorem says that the three intersection
points 
\[ AE \cap BF, \quad BD \cap CE, \quad AD \cap CF \] 
(corresponding to the three minors of the array) are collinear (see
Diagram~\ref{diag:pascal}).
This is one of the oldest theorems in classical projective geometry,
and the picture is so striking that a verbal explanation is almost unnecessary. 

The line containing the intersections is called the Pascal line (or just the pascal)
of the array; we will tentatively denote it as
$\pasc{A}{B}{C}{F}{E}{D}$. The pascal evidently remains the same if we
permute the rows or the columns of the array; thus the same line could
also be written as $\pasc{B}{C}{A}{E}{D}{F}$ or
$\pasc{F}{D}{E}{A}{C}{B}$ etc.

But any essentially different arrangement of the same points,
say $\pasc{B}{F}{D}{C}{A}{E}$,
corresponds \emph{a priori} to a different pascal. 
Since there are $6!$ ways to allocate the points into an array, and $2
\times 3!= 12 $ ways to write
the same pascal, the number of possible pascals is $\frac{6!}{2 \times
  3!} = 60$. It is a theorem due to Pedoe~\cite{Pedoe} that these 
lines are pairwise distinct, if the initial six points are chosen generally.

To recapitulate, six general points on a conic give rise to a
collection of sixty pascals. The entirety of this
configuration\footnote{These sixty pascals 
  satisfy some incidence theorems, leading to the so-called Steiner
  and Kirkman points (see Section~\ref{section.steiner.kirkman.theorems}
  below). These points in turn satisfy some further incidences,
  leading to Cayley lines, Pl{\"u}cker lines and Salmon points. Sometimes these
  auxiliary geometric elements are also included in the term `mystic
  hexagram' (see~\cite{ConwayRyba}). However, they will play no
  substantial role in this paper.} is usually called
Pascal's \emph{hexagrammum mysticum}, or `mystic hexagram'.

\subsection{A dimension count and the enumerative problem}
\label{section.enumerative} 
We will now give an informal description of the enumerative problem to be
considered in this paper. The necessary technical terminology will be
introduced later in Section~\ref{section.pascalvariety}. 

Once and for all, fix the conic $\conic$ inside the complex projective
plane. Recall that $\conic$ has
dimension one (as an algebraic variety), and
hence the set of sextuples $(A, B, \dots, F) \in \conic^6$ has
dimension six. As the sextuple varies in $\conic^6$, so does the
hexagram. 

Now let $\ell$ be a line in $\P^2$, and consider the subset of
sextuples $(A, B, \dots, F)$ such that
\begin{equation} \pasc{A}{B}{C}{F}{E}{D} = \ell.
\label{pascal.constraint} \end{equation} 
Such a subset will later be called a Pascal variety.
Since the lines in $\P^2$ form a two-dimensional family, we expect the
Pascal variety to be of codimension two 
in the set of all sextuples. The same would be true if one substitutes
any of the sixty possible arrays in~(\ref{pascal.constraint}). 

\begin{figure}
\includegraphics[width=12cm]{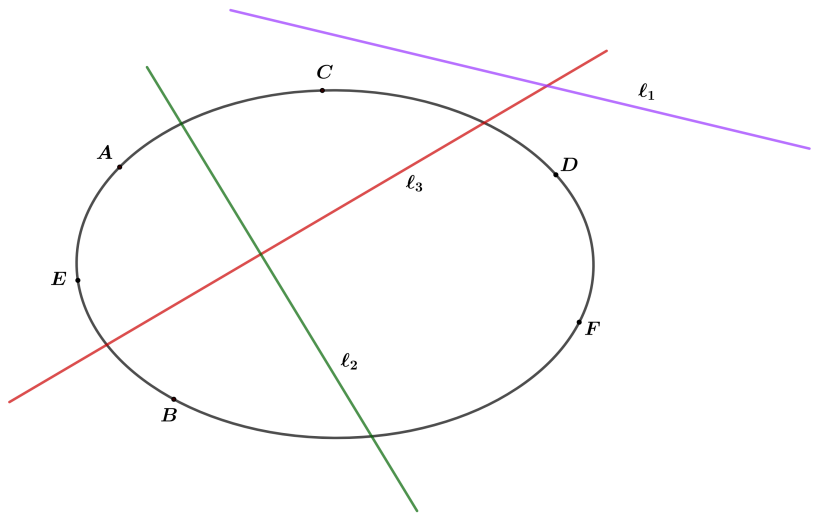}
\caption{Sextuples with pre-specified pascals} 
\label{diag:three_lines} 
\end{figure}

Now fix three general lines $\ell_1, \ell_2, \ell_3$ in the plane, and
(by way of example) consider the subset of sextuples $(A, B, \dots, F)$ such that
\begin{equation}
  \pasc{A}{B}{C}{F}{E}{D} = \ell_1, \qquad
  \pasc{A}{C}{F}{B}{D}{E} = \ell_2, \qquad
  \pasc{D}{A}{F}{B}{C}{E} = \ell_3.
\label{pascal.threeconstraints} \end{equation}
(The three arrays were chosen randomly.)
This subset is of codimension $2 \times 3 = 6$; that is to
say, there is a \emph{finite} number of sextuples which satisfy
the conditions in~(\ref{pascal.threeconstraints}). We should like
to find this number for each such triple of arrays. It will turn out
that there are $77$ essentially distinct triples to be considered. The main result of
this paper is a determination of this number in all of these cases. We
will use computational techniques in commutative algebra to do so. 

The problem is shown schematically in Diagram~\ref{diag:three_lines}.
The three coloured pascals are pre-specified, and we are to find all positions of
$(A,B, \dots, F)$ which make the equalities in
(\ref{pascal.threeconstraints}) come true. There is no obvious
geometric construction which will lead us from the pascals to the
sextuples, and the solution will have to be essentially algebraic. 
\subsection{}
The paper is organised as follows. In the next section, we formulate
the notion of a Pascal variety. The required number is the
cardinality of the intersection of three such varieties.
We then describe the computational procedure to find this
number. In order to list the cases, it is convenient to
introduce the so-called `dual notation' for Pascals, which makes
crucial use of the unique nontrivial outer automorphism of the symmetric group
$\SG_6$. This is explained in Section~\ref{section.labelling}. Using
this notation, we tabulate the results on
page~\pageref{table.77cases} of Section~\ref{section.tripleintersection}.
Some further lines of investigation are
sketched in Section~\ref{section.steiner.kirkman.theorems}. 

All the computations were done using a conjunction of {\sc Maple} and
{\sc Macaulay2}. 

\subsection{References} 
The literature on Pascal's theorem is very large. The standard
classical reference is by George Salmon
(see~\cite[Notes]{SalmonConics}). The dual notation,
together with a host of results discovered by Cremona and
Richmond are explained by H.~F.~Baker in his note `On the 
\emph{Hexagrammum Mysticum} of Pascal' in~\cite[Note II,
pp.~219--236]{Baker}. One of the best modern surveys of this material is by Conway and
Ryba~\cite{ConwayRyba, ConwayRyba2}, which also contains a
bibliography of older literature on this subject. 
We refer the reader to~\cite{Coxeter, KK, Seidenberg} for foundational notions in 
projective geometry; in particular each of them contains a proof of
Pascal's theorem. The modest amount of commutative algebra and
algebraic geometry that we need may be found in~\cite{EisenbudHarris,
  Harris, Smith_etal}.

\subsection{} During the past three centuries, scores of papers on
various aspects of Pascal's hexagram have been published. However, to
the best of my knowledge, this specific enumerative problem has never
been considered before. One possible explanation is that the
computations cannot be done by hand. As we will see below, they
involve calculating the primary decomposition of a rather complicated ideal
inside a polynomial ring in six variables. In my experience, it
usually takes several minutes of computation in {\sc Macaulay2} to work through a single case. For this
reason, several natural questions in this territory still appear to be
inaccessible (see Section~\ref{section.steiner.kirkman.theorems}). 

\section{Pascal varieties}
\label{section.pascalvariety}

In this section, we will introduce the necessary geometric set-up.

\subsection{}
Let $[z_0,z_1,z_2]$ be the homogeneous coordinates in $\P^2$. For instance,
we will denote the line $3 \, z_0  + z_1 + 2
\, z_2 =0$ as $\langle 3, 1, 2 \rangle$.  Identify
$\conic$ with the conic $z_0 \, z_2 = z_1^2$, and fix an isomorphism $\tau:
\P^1 \lra \conic$ by the formula $\tau(r) = [1,r,r^2]$ for $r \in 
\complex$, and $\tau(\infty) = [0, 0, 1]$.

\subsection{} \label{defn.qs}
Consider the set of letters $\ltr = \{\bbA, \bbB, \bbC, \bbD, \bbE,
\bbF \}$. These should be seen formally at the moment, but we
will soon relate them to points on the conic. 
Define a Pascal symbol to be an array such as 
$s = \pasc{\bbE}{\bbC}{\bbA}{\bbB}{\bbF}{\bbD}$, determined up to row and
column shuffles. There are sixty such symbols.

\medskip

Let $\conic^\ltr$ denote the set of maps $\ltr \lra \conic$.
This is a projective variety isomorphic to $(\P^1)^6$.
A \emph{hexad} is an injective map $\ltr \stackrel{h}{\lra}
\conic$; it corresponds to six distinct points
\[ A = h(\bbA), \quad B = h(\bbB), \quad \dots \quad F = h(\bbF) \]
on the conic. If $\Hex$ denotes the set of all hexads, then $\Hex \subseteq
\conic^\ltr$ is a quasiprojective algebraic variety. Given a
Pascal symbol $s$, we get a morphism 
\[ p_s: \Hex \lra \Planedual, \]
which sends a hexad $h$ to the corresponding pascal obtained by
substituting the points $A, \dots, F$ into the symbol $s$. Given
a line $\ell \in \Planedual$, define the Pascal variety 
\begin{equation}
  \Pi(s,\ell) = p_s^{-1}(\{\ell\}) = \{ h \in \Hex: p_s(h) = \ell \}.
  \label{pascal.variety} \end{equation} 

We have the following result about the structure of $\Pi(s,\ell)$.
\begin{Proposition} \rm
The variety $\Pi(s, \ell)$ is isomorphic to a dense open set in
$\conic^4$; in particular it is $4$-dimensional.
\end{Proposition}
\demo 
Without loss of generality, we may assume
\begin{equation} s =
  \pasc{\bbA}{\bbB}{\bbC}{\bbF}{\bbE}{\bbD}. 
  \label{std.symbol}  \end{equation}
First, assume that $\ell$ is not tangent to $\conic$.
Let $\conic \cap \ell = \{P_1, P_2\}$, and write $\conic' = \conic \setminus \{P_1,
P_2\}$. Choose arbitrary distinct points $A, B, C, F$ in $\conic'$,
and let $Q_1 = BF \cap \ell, Q_2 = CF \cap \ell$. Now let $D$
(respectively $E$) be the other intersection of $\conic$ with the
line $AQ_2$ (respectively $AQ_1$). This gives a
hexad in $\Pi(s, \ell)$ and all hexads in it are obtained this way.

If $\ell$ is tangent to $\conic$, then the same argument goes through
with $\conic' = \conic \setminus \{P\}$ where $P$ is the point of
tangency. This proves the result. \qed 

\medskip

Now consider a triple
\begin{equation}
  \label{triple.notation}
  T = \{ (s_1, \ell_1), (s_2, \ell_2), (s_3, \ell_3) \},
\end{equation} 
      where $s_1, s_2, s_3$ are pairwise distinct Pascal symbols and
      $\ell_1, \ell_2, \ell_3$ are general lines in $\P^2$. Our basic
      enumerative problem is to find the number of points in the variety
      \begin{equation} \Pi_T = \bigcap\limits_{i=1}^3 \; \Pi(s_i,
        \ell_i).
        \end{equation}

\subsection{Equations for the Pascal variety}

Fix the symbol $s$ as in~(\ref{std.symbol}). Let $a, b, c, d, e, f$ be
indeterminates, and consider the polynomial ring $R =
\complex[a,b,c,d,e,f]$. Write 
\[ A = \tau(a), \quad B = \tau(b), \quad \dots \quad F = \tau(f). \]
Now it is easy to calculate the coordinates
of the lines $AE, BF$ etc, and eventually those
of the Pascal $\pasc{A}{B}{C}{F}{E}{D}$. They turn out to be 
$\langle u_0, u_1, u_2 \rangle$, where 
\begin{equation} \begin{aligned} 
u_0 & = abde-abdf-acde+acef+bcdf-bcef, \\ 
u_1 & = -abe+abf+acd-acf+adf-aef-bcd+bce-bde+bef+cde-cdf, \\ 
u_2 & = -ad+ae+bd-bf-ce+cf . \end{aligned} \label{coord.pascal} \end{equation} 
(This computation was programmed in {\sc Maple}.)
Given the line $\ell = \langle \alpha_0, \alpha_1, \alpha_2 \rangle$, we have
\[ \pasc{A}{B}{C}{F}{E}{D} = \langle u_0, u_1, u_2 \rangle =
  \langle \alpha_0, \alpha_1, \alpha_2 \rangle, \]
exactly when the $2 \times 2$ minors of the matrix
$\left[ \begin{array}{ccc} u_0 & u_1 & u_2 \\
          \alpha_0 & \alpha_1 & \alpha_2 \end{array} \right]$ are
      zero. These minors generate the ideal $I(s, \ell) \subseteq R$ which locally defines the
      Zariski closure of the Pascal variety $\Pi(s,\ell)$. A
      similar construction works for any symbol $s$. 

      \medskip
      
      Let $\Delta  = \conic^\ltr \setminus \Hex$, which is usually called the big
diagonal. It corresponds to sextuples $(A, \dots, F)$ where some of
the points coincide. By construction, $\Delta$ is disjoint from $\Pi(s,\ell)$.

      \subsection{}      \label{section.computation}
      Given a triple $T$ as above, define the ideal sum 
      \begin{equation} I_T = I(s_1, \ell_1) + I(s_2, \ell_2) + I(s_3,
        \ell_3),
        \label{ideal.sum} \end{equation}
and let $J_T = \sqrt{I_T}$. It is not quite the case that $J_T$ defines $\Pi_T$,
      since some irreducible components of the variety defined by
      it lie inside $\Delta$. In order to remove them, consider the
      minimal prime decomposition 
      \[
        J_T = \bigcap \; \frakp_i ,   \]
      and define
      the restricted intersection $\Lambda_T =\bigcap' \; \frakp_i$, where
      we exclude those
      prime ideals which contain any expression of the form
      $u-v$, for $u,v \in \{a, b, c, d, e, f \}$. If the $\ell_i$
      are chosen generally, then $\Lambda_T \subseteq R$ is an ideal of
      height $6$. Now we have\footnote{The only points of
        $\conic^\ltr$ which lie outside the open subset $\text{Spec}
        \, R$ are those where one of the coordinates is
        $\infty$. If the $\ell_i$ are chosen generally, none of them
        will be included in $\Pi_T$.} 
      \begin{equation} \text{card} \, (\Pi_T) = \deg \Lambda_T.
      \label{deg.lambdat} \end{equation} 
      We denote this number by $\IN_T$, since it is independent of the
      $\ell_i$ as long as the lines are chosen generally. With this
      understood, we will simply write $T = \{s_1, s_2, s_3\}$ instead of~(\ref{triple.notation}).  

      \subsection{} I have used {\sc Maple} to program the ideal
      sum in (\ref{ideal.sum}), and then programmed the rest of the
      computation in {\sc Macaulay2}. There are several techniques
      for calculating the radical and the primary decomposition of an
      ideal inside a polynomial ring. By and large, they depend crucially
      upon the concept of
      a Gr{\"o}bner basis. The reader will find an explanation and 
      comparison of such techniques in~\cite{DGP}.

      There are practical difficulties arising from
      the fact that $I_T$ is generated by $9$ large polynomials in the variables
      $a, \dots, f$. This leads to a corresponding difficulty in calculating the prime
      decomposition of its radical. Hence, this is the procedure I have followed: instead of using
      ${\mathbb Q}$ as the base field (which slows down the computation)
      I have worked through each case thrice by using large
      prime fields $\FF_p$ (such as $p = 32003, 43051, 48619$) and varying the
      choice of $\ell_i$ each time. Since the results are
      consistent, I have a high degree of confidence in their
      correctness.

      In order to tabulate the results in a compact form, we will use
      the dual notation for pascals. It is important
      in its own right, since it is used throughout all later
      development of the subject. Be that as it may, the mathematics 
      underlying it has a charm of its own. 
            
      \section{The dual notation for Pascals}
      \label{section.labelling} 

The symmetric group $\SG_6$ has a unique outer automorphism, which
turns out to have a close connection to the geometry of Pascal's
hexagram. It is the basis of the labelling scheme explained below
(which is a minor variant of the one in Baker~\cite{Baker}).
The notation introduced by Conway
and Ryba~\cite{ConwayRyba2} is ostensibly different, but based upon the
same principle.

\subsection{} 
Recall that if $G$ is a group and $g \in G$, then 
\[ G \lra G, \qquad x \lra g^{-1} \, x \, g \] 
is called an inner automorphism of $G$. The subgroup of inner
automorphisms (denoted $\Inn(G)$) is a normal subgroup of $\Aut(G)$
(the group of all automorphisms), and the quotient 
$\Out(G) = \frac{\Aut(G)}{\Inn(G)}$ 
is the group of outer automorphisms. The automorphisms of the
symmetric groups $\SG_n$ are characterised by the following theorem due to
H{\"o}lder (see \cite[Ch.~7]{Rotman}):  
\[ \Out(\SG_n) = 
\begin{cases} 
\{ e\} & \text{for $n \neq 6$,} \\ 
\; \ZZ_2 & \text{for $n =6$.} 
\end{cases} 
\] 
Thus, $\SG_6$ has an essentially unique nontrivial outer automorphism. We will not
reproduce its construction here, since it may be found in several
places (e.g., see~\cite{HMSV, JanuszRotman}). For our purposes,
it will be more convenient to see it 
as an isomorphism between two different copies of $\SG_6$. 
\subsection{} 
\label{section.outerauto} 
In addition to $\ltr = \{\bbA,\bbB,\bbC,\bbD,\bbE,\bbF\}$,
consider the set $\six = \{1,2,3,4,5,6\}$. 
For any set $X$, let $\SG(X)$ denote the group of bijections $X
\ra X$. Then the following table gives an isomorphism
\[ \zeta: \SG(\ltr) \lra \SG(\six), \]
which is one realisation of the outer automorphism. For instance, it is to be understood as saying that 
$\zeta$ takes the transposition $\transp{\bbA}{\bbB}$ to the cycle 
$\cyclett{1}{4}{2}{5}{3}{6}$ of type $2 + 2 + 2$. 
\[ 
\begin{array}{|c|c|c|c|c|c|c|} \hline 
{} & \bbB & \bbC & \bbD & \bbE & \bbF \\ \hline 
\bbA & 14.25.36 & 16.24.35 & 13.26.45 & 12.34.56 & 15.23.46 \\ 
\bbB & {} & 15.26.34 & 12.35.46 & 16.23.45 & 13.24.56 \\ 
\bbC & {} & {} & 14.23.56 & 13.25.46 & 12.36.45 \\ 
\bbD & {} & {} & {} & 15.24.36 & 16.25.34 \\ 
\bbE & {} & {} & {} & {} & 14.26.35 \\ \hline 
\end{array} \label{table.zeta} \]  

\medskip 

Similarly, $\zeta^{-1}$ is given by the following table: 
\[ 
\begin{array}{|c|c|c|c|c|c|c|} \hline 
{} & 2 & 3 & 4 & 5 & 6 \\ \hline 
1 & \synth{\bbA}{\bbE}{\bbB}{\bbD}{\bbC}{\bbF} & \synth{\bbA}{\bbD}{\bbB}{\bbF}{\bbC}{\bbE}  & 
\synth{\bbA}{\bbB}{\bbC}{\bbD}{\bbE}{\bbF}  & \synth{\bbA}{\bbF}{\bbB}{\bbC}{\bbD}{\bbE} & \synth{\bbA}{\bbC}{\bbB}{\bbE}{\bbD}{\bbF} \\ 
2 & {} & \synth{\bbA}{\bbF}{\bbB}{\bbE}{\bbC}{\bbD} & 
\synth{\bbA}{\bbC}{\bbB}{\bbF}{\bbD}{\bbE} & \synth{\bbA}{\bbB}{\bbC}{\bbE}{\bbD}{\bbF}  & \synth{\bbA}{\bbD}{\bbB}{\bbC}{\bbE}{\bbF}  \\ 
3 & {} & {} & \synth{\bbA}{\bbE}{\bbB}{\bbC}{\bbD}{\bbF} &
\synth{\bbA}{\bbC}{\bbB}{\bbD}{\bbE}{\bbF} & \synth{\bbA}{\bbB}{\bbC}{\bbF}{\bbD}{\bbE} \\ 
4 & {} & {} & {} &  \synth{\bbA}{\bbD}{\bbB}{\bbE}{\bbC}{\bbF} & \synth{\bbA}{\bbF}{\bbB}{\bbD}{\bbC}{\bbE} \\ 
5 & {} & {} & {} &  {} & \synth{\bbA}{\bbE}{\bbB}{\bbF}{\bbC}{\bbD} \\ \hline 
\end{array} \]

\smallskip

For instance, it takes $\transp{1}{2}$ to
$\cyclett{\bbA}{\bbE}{\bbB}{\bbD}{\bbC}{\bbF}$.

\subsection{} \label{section.pascal.labelling}
Our labelling schema uses the fact that 
the $15$ elements of type $2+2+2$ in $\SG(\ltr)$ are bijectively 
mapped by $\zeta$ onto the $15$ transpositions in $\SG(\six)$.

By way of example, consider the Pascal symbol
$\pasc{\bbA}{\bbD}{\bbE}{\bbC}{\bbB}{\bbF}$. There are six
lines going between the rows; separate them into two diagrams
as follows:

\[ \begin{aligned} 
\xymatrix{ 
\bbA \ar@{-}[dr] & \bbD \ar@{-}[dr] & \bbE \ar@{-}[dll]\\ 
\bbC & \bbB & \bbF } & \qquad & \xymatrix{ 
\bbA \ar@{-}[drr] & \bbD \ar@{-}[dl] & \bbE \ar@{-}[dl]\\ 
\bbC & \bbB & \bbF } \\ \\
\end{aligned} \] 
In either of the diagrams, each letter has a unique edge 
joining it to another letter lying in a different column. Read both
diagrams as cycles of type $2+2+2$, and observe that 
\[ \cyclett{\bbA}{\bbB}{\bbD}{\bbF}{\bbC}{\bbE} \stackrel{\zeta}{\lra} (2 \, 5), \qquad 
\cyclett{\bbA}{\bbF}{\bbD}{\bbC}{\bbE}{\bbB} \stackrel{\zeta}{\lra} (2 \, 3). \] 
Single out the $2$ common to both cycles, leaving
$3$ and $5$. Hence we label this pascal as $k(2,\{3,5\})$. It is
customary to flatten this out for readability, and 
write it indifferently as either $k(2,35)$ or $k(2,53)$. A moment's
reflection will show that a row or column shuffle in a symbol does not change its
label. 

\medskip 

Thus any pascal is labelled as $k(w,xy)$ for some $w \in \six$, and $x,y \in \six \setminus
\{w\}$, where the order of $x,y$ is immaterial.
There are $6 \times \binom{5}{2} = 60$ such labels, bijectively
corresponding to the Pascal symbols. It is easy to reconstruct the array from the
label; for instance, starting from $k(3,15) = k(3,51) $, 
\[ (3 \, 1) \stackrel{\zeta^{-1}}{\lra} \cyclett{\bbA}{\bbD}{\bbB}{\bbF}{\bbC}{\bbE}, 
\qquad (3 \, 5) \stackrel{\zeta^{-1}}{\lra} \cyclett{\bbA}{\bbC}{\bbB}{\bbD}{\bbE}{\bbF}, \] 
which must correspond to $\pasc{\bbA}{\bbB}{\bbE}{\bbF}{\bbC}{\bbD}$. Of course, the
array is determined only up to row and column shuffles. As an
exercise, the reader may
wish to check that $\pasc{\bbA}{\bbB}{\bbC}{\bbF}{\bbE}{\bbD}$
corresponds to $k(1,23)$.

\section{The triple intersections}
\label{section.tripleintersection}

\subsection{} The group $\SG(\ltr)$ acts on $\conic^\ltr$ and hence on
$\Hex$ in an obvious way. Let $\Omega$ denote the set of unordered triples
      of Pascal symbols. This set has $\binom{60}{3} = 34,220$
      elements. For instance,
      \[ T = \left\{ \pasc{\bbA}{\bbB}{\bbC}{\bbF}{\bbE}{\bbD}, 
          \pasc{\bbA}{\bbC}{\bbF}{\bbB}{\bbD}{\bbE}, 
          \pasc{\bbD}{\bbA}{\bbF}{\bbB}{\bbC}{\bbE} \right\} \]
      is an element of $\Omega$. It is clear that $\SG(\ltr)$ acts on
      $\Omega$ as well. If $\sigma \in \SG(\ltr)$ takes the triple
      $T_1$ to $T_2$, then we have a bijection
      \[ \Pi_{T_1} \stackrel{\sim}{\lra} \Pi_{T_2}, \qquad h \lra
        h^\sigma \]
      for a fixed choice of $\ell_1, \ell_2, \ell_3$. 
      In other words, every hexad in $\Pi_{T_2}$ is obtained by
      relabelling a unique hexad in $\Pi_{T_1}$ as dictated by the
      permutation $\sigma$.  It follows that $\IN_{T_1} =
      \IN_{T_2}$. Hence it is only necessary to decompose $\Omega$ into its
      $\SG(\ltr)$-orbits, and find $\IN_T$ for any
      one triple $T$ in each orbit. 

      Using the labelling scheme, we might as well regard $\Omega$ as the
      set of unordered triples of Pascal labels with a natural action
      of $\SG(\six)$. For instance, the triple written above now appears
      in a more compact form as
      \[ \left\{  \, k(1,23), \; k(3,14), \; k(6, 35) \right\}. \]
      It is routine to program
      this group action in {\sc Maple}. It turns out that there are
      altogether $77$ orbits. In each case, I have carried out the 
      computational procedure from Section~\ref{section.computation}
      on lexicographically\footnote{We describe this order
        in brief, although its details are of minor importance. 
        Within each triple, elements of the form $k(1,xy)$ are
        listed first, followed by those of the form $k(2,xy)$ and so
        on. Ties are broken by writing $x, y$ in increasing order. The
        order between two triples is decided by comparing the first
        element in each, following by the second in each, and so on.}
      the smallest element in each orbit.

      \medskip
      
The results are summarised in the following theorem, which refers to
the table on page~\pageref{table.77cases}.

\begin{Theorem} \rm 
With notation as above: 
\begin{itemize} 
\item 
Any triple of Pascal labels lies in the
$\SG(\six)$-orbit of exactly one triple appearing\footnote{The redundant '$k$' is
omitted throughout.} in the table. 
\item 
For any such triple $T$, the intersection number $\IN_T$ is as shown in the table. 
\end{itemize} 
\end{Theorem}

The triples are in lexicographic order when read from left to
right and top to bottom. The first label is always
$k(1,23)$. The colour-coding of the entries has the
following meaning: 
\begin{itemize}
\item
  The three \textcolor{red}{red} entries have intersection number zero.  The geometry
  of these cases is discussed in
  Section~\ref{section.steiner.kirkman.theorems} below.
\item \label{example.case}
  If $T = \{s_1, s_2, s_3\}$ is a triple, then define its stabiliser
  $G_T$ to be the set of elements $\sigma \in
  \SG(\six)$ such that $\sigma(s_i) = s_i$ for $i=1,2,3$. It is easy to
  determine this group for each triple; then it turns out that the only
  possible nontrivial stabilisers are $\ZZ_2$ and $\ZZ_2
  \times \ZZ_2$. The \textcolor{blue}{blue} entries are the ones with $G_T
  \simeq {\mathbb Z}_2$. For instance, if 
  \[ T = \{(1,23), (1, 45), (2,45)\}  \]
  (the entry in row 10, column 2), then $G_T = \{e, (4 \, 5)
  \}$. A complete example for this case is
  given in Section~\ref{section.example} below. Notice that 
  \[ (4 \, 5) \stackrel{\zeta^{-1}}{\lra}
    (\bbA \, \bbD)   (\bbB \, \bbE)   (\bbC \, \bbF). 
  \]
 Hence, starting from a hexad in $\Pi_T$, we get another one by
  the simultaneous interchange
  \[ a \leftrightarrow d, \quad b \leftrightarrow e, \quad c
    \leftrightarrow f. \] 
  This implies that the $8$-element set $\Pi_T$ can be partitioned into
  four subsets of two elements each, where the hexads in each subset
  are obtained from each other by this interchange. 
  Of course, a similar remark applies to each \textcolor{blue}{blue} entry. 
\item 
  There is a unique triple $Q$ with stabiliser $\ZZ_2 \times \ZZ_2$. It
  appears in row 11, column 2 and is coloured \textcolor{brown}{brown}. It is easy to
  see that
  \[ G_Q =  \{e, (2 \, 3), (4 \, 5), (2 \, 3) (4 \, 5) \}. \]
 As above, this implies that the $8$-element set $\Pi_Q$ can be 
partitioned into two subsets of four elements each, where the hexads in
each subset are permuted by the elements of the
stabiliser. 
\end{itemize}

\smallskip

The entry in row 25, column 3 is especially simple and elegant. It
  corresponds to the smallest nonzero intersection number. This case
is treated at length in \cite{AC}. Its discovery was something of an  
inspiration for this paper. 

\begin{table} \caption*{\label{table.77cases}
    {\sc{The table of triple-intersection numbers}}} 
  \begin{tabular}{||r|| cr|cr|cr|} \hline \hline
1 & (1, 23), (1, 24), (1, 25) &  22 & 
\textcolor{red}{(1, 23), (1, 24), (1, 34)} & \textcolor{red}{0} & 
(1, 23), (1, 24), (1, 35) &  20 \\ 
2 & \textcolor{blue}{(1, 23), (1, 24), (1, 56)} &  \textcolor{blue}{12} & 
\textcolor{blue}{(1, 23), (1, 24), (2, 13)} &   \textcolor{blue}{18} & 
(1, 23), (1, 24), (2, 15) &  8 \\ 
3 & \textcolor{blue}{(1, 23), (1, 24), (2, 34)} & \textcolor{blue}{22} & 
(1, 23), (1, 24), (2, 35) & 28 & 
\textcolor{blue}{(1, 23), (1, 24), (2, 56)} &  \textcolor{blue}{22} \\ 
4 & \textcolor{blue}{(1, 23), (1, 24), (3, 12)} &  \textcolor{blue}{28} & 
\textcolor{blue}{(1, 23), (1, 24), (3, 14)} &   \textcolor{blue}{20} & 
(1, 23), (1, 24), (3, 15) &  12 \\ \hline
5 & \textcolor{blue}{(1, 23), (1, 24), (3, 24)} &   \textcolor{blue}{16} & 
(1, 23), (1, 24), (3, 25) &  30 & 
(1, 23), (1, 24), (3, 45) &   28 \\ 
6 & \textcolor{blue}{(1, 23), (1, 24), (3, 56)} & \textcolor{blue}{16} & 
(1, 23), (1, 24), (5, 12) &  30 & 
(1, 23), (1, 24), (5, 13) & 28 \\ 
7 & (1, 23), (1, 24), (5, 16) &  14 & 
(1, 23), (1, 24), (5, 23) &  14 & 
(1, 23), (1, 24), (5, 26) & 22 \\ 
8 & (1, 23), (1, 24), (5, 34)& 18 & 
(1, 23), (1, 24), (5, 36) & 20 & 
\textcolor{blue}{(1, 23), (1, 45), (2, 13)} &  \textcolor{blue}{16} \\ \hline
9 & (1, 23), (1, 45), (2, 14) &  12 & 
\textcolor{blue}{(1, 23), (1, 45), (2, 16)} &  \textcolor{blue}{12} & 
(1, 23), (1, 45), (2, 34)&  20 \\ 
10 & \textcolor{blue}{(1, 23), (1, 45), (2, 36)} &  \textcolor{blue}{12} & 
\textcolor{blue}{(1, 23), (1, 45), (2, 45)} &  \textcolor{blue}{8} & 
(1, 23), (1, 45), (2, 46)&  18 \\ 
11 & \textcolor{blue}{(1, 23), (1, 45), (6, 12)} &  \textcolor{blue}{16} & 
\textcolor{brown}{(1, 23), (1, 45), (6, 23)} &  \textcolor{brown}{8} & 
(1, 23), (1, 45), (6, 24)&  20 \\ 
12 & \textcolor{red}{(1, 23), (2, 13), (3, 12)} &  \textcolor{red}{0} & 
\textcolor{blue}{(1, 23), (2, 13), (3, 14)} &  \textcolor{blue}{18} & 
\textcolor{blue}{(1, 23), (2, 13), (3, 45)} &  \textcolor{blue}{16} \\ \hline 
13 & \textcolor{blue}{(1, 23), (2, 13), (4, 12)} &  \textcolor{blue}{20} & 
\textcolor{blue}{(1, 23), (2, 13), (4, 13)} &  \textcolor{blue}{14} & 
(1, 23), (2, 13), (4, 15)&  26 \\ 
14 & (1, 23), (2, 13), (4, 35)&  28 & 
\textcolor{blue}{(1, 23), (2, 13), (4, 56)} &  \textcolor{blue}{8} & 
\textcolor{red}{(1, 23), (2, 14), (3, 14)} &  \textcolor{red}{0} \\ 
15 & (1, 23), (2, 14), (3, 15)&  8 & 
\textcolor{blue}{(1, 23), (2, 14), (3, 24)} & \textcolor{blue}{20} &  
(1, 23), (2, 14), (3, 25)&  14 \\ 
16 & (1, 23), (2, 14), (3, 45)& 12 &  
\textcolor{blue}{(1, 23), (2, 14), (3, 56)} &  \textcolor{blue}{12} & 
(1, 23), (2, 14), (5, 12)& 20 \\ \hline
17 & (1, 23), (2, 14), (5, 13)&  20 & 
(1, 23), (2, 14), (5, 14)&  4 & 
(1, 23), (2, 14), (5, 16)& 16 \\ 
18 & (1, 23), (2, 14), (5, 34)&  20 & 
(1, 23), (2, 14), (5, 36)&  12 & 
(1, 23), (2, 34), (3, 45)&  28 \\ 
19 & \textcolor{blue}{(1, 23), (2, 34), (3, 56)} & \textcolor{blue}{20} & 
\textcolor{blue}{(1, 23), (2, 34), (4, 13)} &  \textcolor{blue}{30} & 
(1, 23), (2, 34), (4, 15)& 22 \\ 
20 & (1, 23), (2, 34), (4, 35)& 30 & 
\textcolor{blue}{(1, 23), (2, 34), (4, 56)} & \textcolor{blue}{16} & 
(1, 23), (2, 34), (5, 14)& 28 \\ \hline
21 & (1, 23), (2, 34), (5, 16)& 14 & 
(1, 23), (2, 34), (5, 23)&  12 & 
(1, 23), (2, 34), (5, 24)& 24 \\ 
22 & (1, 23), (2, 34), (5, 26)& 28 & 
(1, 23), (2, 34), (5, 34)&  16 & 
(1, 23), (2, 34), (5, 36)&  30 \\ 
23 & (1, 23), (2, 34), (5, 46)&  20 & 
\textcolor{blue}{(1, 23), (2, 45), (3, 45)} & \textcolor{blue}{6} & 
(1, 23), (2, 45), (3, 46)& 12 \\ 
24 & (1, 23), (2, 45), (4, 16)&  24 & 
(1, 23), (2, 45), (4, 36)& 24 & 
\textcolor{blue}{(1, 23), (2, 45), (6, 23)} & \textcolor{blue}{8} \\
     \hline
25 & (1, 23), (2, 45), (6, 34)&  28 & 
\textcolor{blue}{(1, 23), (2, 45), (6, 45)} &  \textcolor{blue}{8} & 
\textcolor{blue}{(1, 23), (4, 23), (5, 23)} & \textcolor{blue}{2} \\ 
26 & (1, 23), (4, 23), (5, 26)&  14 & 
                                 (1, 23), (4, 25), (6, 35) &  24 & {} &
     \\ \hline 
                                                                  
  \end{tabular} \end{table}

\newpage

\subsection{} 
\label{section.example} 
Given a triple $T = \{(s_1, \ell_1), (s_2, \ell_2), (s_3, \ell_3)
\}$, one should like to write down the set $\Pi_T$
explicitly at least as an illustration. This is awkward to do
over the complex numbers, since
the generators of the ideal $\Lambda_T$ are usually unwieldy. As
a compromise, I have written down such an example where the base field
is $\overline{\FF}_{101}$, namely the algebraic closure of the
finite field with $101$ elements. The prime is large enough to give an
interesting example, and small enough to keep the coefficients tidy. 

Let
\[ s_1 = (1, 23), \quad s_2 = (1,45), \quad s_3 = (2,45), \]
and
\[ \ell_1 = \langle 1, 35, 48 \rangle, \quad
  \ell_2 = \langle 1, 5, 26 \rangle, \quad
  \ell_3 = \langle 1, 32, 52 \rangle. \]

As explained on page~\pageref{example.case}, the set $\Pi_T$
breaks up into four $G_T$-orbits of two elements each. Within each orbit, we get one
solution from another by simultaneously
interchanging $a$ with $d$, $b$ with $e$, and $c$
with $f$. I have written down the elements in $\Pi_T$ by following the procedure in
Section~\ref{section.computation}. 

The first orbit is given by the solution 
\begin{equation}
  a = 48, \quad b = 49, \quad c = 14, \quad d = 92, \quad e = 9,
  \quad f = 57, 
  \label{orbit1} \end{equation}
together with its interchange. 

The second orbit is given by 
\begin{equation}
    b = 29 \, a + 69, \quad c = 58 \, a + 70, \quad d = 4-a, \quad e =
    72 \, a + 84, \quad f = 43 \, a + 100,
\label{orbit2}   \end{equation} 
where $a$ is either of the roots of the equation $x^2 - 4 \, x + 63
=0$. Observe that the two roots add up to $4$. If we
substitute $4-a$ for $a$ throughout, then $b, e$ get exchanged and so
do $c, f$. 

The third orbit is given by 
\begin{equation}
    b = 31 \, a + 64, \quad c = 9 \, a + 70, \quad d = 51 - a, \quad e
  = 70 \, a + 29, \quad f = 92 \, a + 24,
\label{orbit3}   \end{equation} 
where $a$ is either of the roots of the equation $x^2 - 51 \, x +
4=0$. The exchange property is exactly as above when $a$ is replaced
by $51-a$. 

The fourth orbit is given by 
\begin{equation}
    b = 43 \, a + 18, \quad c = 71 \, a + 8, \quad d = 56 - a, \quad e
    = 58 \, a + 2, \quad f = 30 \, a + 45,
  \label{orbit4} \end{equation} 
where $a$ is either of the roots of the equation $x^2 - 56 \, x +
4=0$. The exchange property is exactly as above when $a$ is replaced
by $56-a$.

This pattern of solutions is not an intrinsic feature of the triple,
but depends on the choice of the prime and the coefficients of the
$\ell_i$. For instance, the polynomial
$x^2 - 4 \, x + 63$ in the second orbit does not split over $\FF_{101}$ because its
discriminant $-236$ is a non-residue modulo $101$. If it did, the
second orbit would have looked like the first. One expects to see 
different patterns in different examples; it is only the cardinality
of $\Pi_T$ which is intrinsic to $T$. 

\section{Theorems of Steiner and Kirkman}
\label{section.steiner.kirkman.theorems} 

\subsection{Steiner's theorem}
It is a theorem due to Steiner that the pascals
\[ s_1 = k(1,23), \quad s_2 = k(2,13), \quad s_3 = k(3,12), \]
are concurrent (see~\cite{ConwayRyba}). Of course, this is understood to apply to any three
pascals following the same pattern; that is to say,
\[ k(x,yz), \quad k(y,xz), \quad k(z,xy) \]
are concurrent for any $x, y, z \in \six$. This explains why the entry in row 12,
column 1 is zero, since three general lines are not concurrent.

However, this suggests the following modification of the original
set-up: choose a general point $P$ in the plane and three lines
$\ell_1, \ell_2, \ell_3$ passing through $P$. Now consider the variety
\[ X = \bigcap\limits_{i=1}^3 \; \Pi(s_i, \ell_i). \]
A simple count shows that $X$ is of codimension five in $\Hex$; that
is to say, it is a curve. (This is so because $P$ has two degrees of freedom and each
$\ell_i$ has one, leading to a total of $5$.) Since $X \subseteq \conic^\ltr$, there is a natural
projection map $X \stackrel{\pi_L}{\lra} \conic$ for each of the
letters $L \in \ltr$.

I have calculated the equations of $X$ in a few cases. They show
that $X \stackrel{\pi_L}{\lra} \conic$ is a
$7$ to $1$ map for any $L$. In other words, if any one letter (say
$a$) is given a
fixed value, then there are exactly $7$ such points in $X$.

\subsection{Kirkman's theorem}
It is a theorem due to Kirkman that the pascals
\[ s_1 = k(1,23), \quad s_2 = k(1,24), \quad s_3 = k(1,34), \]
are concurrent (see~\cite{ConwayRyba}). As before, this is understood to apply to any three
pascals following the same pattern; that is to say,
\[ k(x,yz), \quad k(x,yw), \quad k(x, zw) \]
are concurrent for any $x, y, z, w \in \six$. This explains
why the entry in row 1, column 2 is zero. 

We can construct a curve 
\[ Y = \bigcap\limits_{i=1}^3 \; \Pi(s_i, \ell_i). \]
as before. Calculations show that 
$Y \stackrel{\pi_L}{\lra} \conic$ is a $4$ to $1$ map for any letter
$L$.

One should like to know more about the geometry of the curves $X$ and
$Y$, and specifically whether they are rational. However, their
equations are far too complicated to make this feasible at the
moment. 

\subsection{}
The triple in the $14$th row and $3$rd column 
corresponds to the arrays
\[ \pasc{A}{B}{C}{F}{E}{D}, \qquad
  \pasc{A}{D}{F}{B}{C}{E}, \qquad
  \pasc{A}{C}{F}{B}{D}{E}. \]
These pascals are concurrent, since they
all pass through $AE \cap BF$. Hence the corresponding entry is
zero. It is possible to construct a curve as in the previous two
sections; we leave these details to the reader. 

\medskip

In the end, we will outline some further lines of investigation. 

\subsection{Steiner and Kirkman varieties} The common point which lies on 
\[ k(x,yz), \quad k(y, xz), \quad k(z,xy), \]
is called the Steiner point
corresponding the set $\{x, y, z \} \subseteq \six$. It is customary
to write it simply as $G[xyz]$, where the order of $x,y,z$ is
immaterial. There are $\binom{6}{3} = 20$ such points.

Given $\{x,y,z\}$ and a general point $P$, it is natural to consider the
Steiner variety consisting of sextuples $(A,B, \dots, F)$ such that
$G[xyz] = P$. This is expected to be of codimension two, and as before
we have the enumerative question of finding the intersection degree of
three such varieties.

Similarly, the point common to the pascals
\[ k(x,yz), \quad k(x,yw), \quad k(x, zw) \] is called the Kirkman point
$K[x,yzw]$. There are $60$ such points, and exactly the same question
can be formulated for them.

\subsection{}
There is another (somewhat speculative) approach to the enumerative
problem considered in this paper. The relevant geometric terminology
can be found in~\cite{EisenbudHarris2, Hartshorne}.

The pascal
$\pasc{A}{B}{C}{F}{E}{D}$ still remains well-defined if the sextuple
$(A, B, \dots, F)$ has a double point, i.e., if two of the points come
together on the conic. However, it may no longer be well-defined
if we have a triple point or two double points. (A precise statement of
possibilities is given in~\cite[\S 2]{CD}.) It follows that
\[ p_s: \conic^\ltr \, - \rightarrow \Planedual \]
from Section~\ref{defn.qs} is only a rational map, and its
indeterminacy locus consists of a union of polydiagonals. It may be possible to use
successive blow-ups (cf.~\cite[Ch.~II.7]{Hartshorne})
in order to construct a smooth projective variety $\Gamma$ such that,
\[ p_s: \Gamma  \longrightarrow \Planedual \]
is a morphism, and $\Gamma$ contains $\Hex$ as a dense open subset on
    which $p_s$ agrees with the earlier geometric
    definition.\footnote{In the paper~\cite{CD},
      Sergio Da Silva and I have partly succeeded in resolving the
indeterminacies of the rational map above. However, for reasons
outlined there, this falls short of a full solution.}
Once this is accomplished, $p_s^{-1}(\{\ell\})$ would be a codimension two
algebraic cycle on
$\Gamma$. Then, assuming that the Chow ring of $\Gamma$ can be
effectively computed, one can interpret $\IN_T$ as the
degree of the intersection of three such cycles. This approach, if
carried out in full, would prove to be interesting and
informative. 

\medskip

As mentioned earlier, this is very likely the
first paper on this enumerative problem. Much remains to be done.

\medskip 

\centerline{--} 

\end{document}